\documentclass[12pt]{amsart}
\def\Image{\operatorname{Im}}
\def\N{{\Bbb N}}

\newtheorem{Theorem}{Theorem}[section]

\newtheorem{Lemma}[Theorem]{Lemma}

\theoremstyle{definition}
\newtheorem{Definition}[Theorem]{Definition}

\theoremstyle{remark}

\begin{document}
\sloppy
\title{On dense orbits in the boundary of a Coxeter system}
\author{Tetsuya Hosaka} 
\address{Department of Mathematics, Utsunomiya University, 
Utsunomiya, 321-8505, Japan}
\date{February 1, 2005}
\email{hosaka@cc.utsunomiya-u.ac.jp}
\thanks{
Partly supported by the Grant-in-Aid for Young Scientists (B), 
The Ministry of Education, Culture, Sports, Science and Technology, Japan.
(No.~15740029).}
\maketitle
\begin{abstract}
In this paper, 
we study the minimality of the boundary of a Coxeter system.
We show that for a Coxeter system $(W,S)$ 
if there exist a maximal spherical subset $T$ of $S$ 
and an element $s_0\in S$ such that $m(s_0,t)\ge 3$ for each $t\in T$
and $m(s_0,t_0)=\infty$ for some $t_0\in T$, 
then every orbit $W\alpha$ is dense 
in the boundary $\partial\Sigma(W,S)$ of the Coxeter system $(W,S)$, 
hence $\partial\Sigma(W,S)$ is minimal, 
where $m(s_0,t)$ is the order of $s_0t$ in $W$.
\end{abstract}

\section{Introduction and preliminaries}

The purpose of this paper is to study 
the minimality of the boundary of a Coxeter system.
A {\it Coxeter group} is a group $W$ having a presentation
$$\langle \,S \, | \, (st)^{m(s,t)}=1 \ \text{for}\ s,t \in S \,
\rangle,$$ 
where $S$ is a finite set and 
$m:S \times S \rightarrow \N \cup \{\infty\}$ is a function 
satisfying the following conditions:
\begin{enumerate}
\item[(1)] $m(s,t)=m(t,s)$ for each $s,t \in S$,
\item[(2)] $m(s,s)=1$ for each $s \in S$, and
\item[(3)] $m(s,t) \ge 2$ for each $s,t \in S$
such that $s\neq t$.
\end{enumerate}
The pair $(W,S)$ is called a {\it Coxeter system}.
Let $(W,S)$ be a Coxeter system.
For a subset $T \subset S$, 
$W_T$ is defined as the subgroup of $W$ generated by $T$, 
and called a {\it parabolic subgroup}.
If $T$ is the empty set, then $W_T$ is the trivial group.
A subset $T\subset S$ is called a {\it spherical subset of $S$}, 
if the parabolic subgroup $W_T$ is finite.

Every Coxeter system $(W,S)$ determines 
a {\it Davis-Moussong complex} $\Sigma(W,S)$ 
which is a CAT(0) geodesic space (\cite{D1}, \cite{D2}, \cite{D3}, \cite{M}).
Here the $1$-skeleton of $\Sigma(W,S)$ is 
the Cayley graph of $W$ with respect to $S$.
The natural action of $W$ on $\Sigma(W,S)$ is proper, cocompact and by isometry.
We can consider a certain fundamental domain $K(W,S)$ 
which is called a {\it chamber} of $\Sigma(W,S)$ 
such that $W K(W,S)=\Sigma(W,S)$ (\cite{D2}, \cite{D3}).
If $W$ is infinite, then $\Sigma(W,S)$ is noncompact and 
$\Sigma(W,S)$ can be compactified by adding its ideal boundary
$\partial \Sigma(W,S)$ (\cite{BH}, \cite[\S 4]{D2}).
This boundary 
$\partial \Sigma(W,S)$ is called the {\it boundary of} $(W,S)$.
We note that the natural action of $W$ on $\Sigma(W,S)$ 
induces an action of $W$ on $\partial \Sigma(W,S)$.

The following theorem was proved in \cite{Ho2}.

\begin{Theorem}\label{Thm1}
Let $(W,S)$ be a Coxeter system. 
Suppose that there exist a maximal spherical subset $T$ of $S$ 
and an element $s_0\in S$ such that $m(s_0,t)\ge 3$ for each $t\in T$
and $m(s_0,t_0)=\infty$ for some $t_0\in T$. 
Then $W\alpha$ is dense in $\partial\Sigma(W,S)$ 
for some $\alpha \in \partial\Sigma(W,S)$.
\end{Theorem}

Suppose that 
a group $G$ acts on a compact metric space $X$ by homeomorphisms.
Then $X$ is said to be {\it minimal}, 
if every orbit $Gx$ is dense in $X$.

For a negatively curved group $\Gamma$ and 
the boundary $\partial\Gamma$ of $\Gamma$, by an easy argument, 
we can show that $\Gamma\alpha$ is dense in $\partial\Gamma$ 
for each $\alpha\in\partial\Gamma$, 
that is, $\partial\Gamma$ is minimal.

In this paper, 
we prove the following theorem as an extension of Theorem~\ref{Thm1}.

\begin{Theorem}\label{Thm}
Let $(W,S)$ be a Coxeter system 
which satisfies the condition in Theorem~\ref{Thm1}.
Then every orbit $W\alpha$ is dense in $\partial\Sigma(W,S)$, 
that is, $\partial\Sigma(W,S)$ is minimal.
\end{Theorem}

\section{Lemmas on Coxeter groups and the Davis-Moussong complexes}

In this section, we recall and prove some lemmas for Coxeter groups 
and the Davis-Moussong complexes which are used later.

\begin{Definition}
Let $(W,S)$ be a Coxeter system and $w\in W$.
A representation $w=s_1\cdots s_l$ ($s_i \in S$) is said to be 
{\it reduced}, if $\ell(w)=l$, 
where $\ell(w)$ is the minimum length of 
word in $S$ which represents $w$.
\end{Definition}

The following lemma is known.

\begin{Lemma}[\cite{Bo}, \cite{Br}, \cite{D1}, \cite{D3}, \cite{Hu}]\label{lem1}
Let $(W,S)$ be a Coxeter system.
\begin{enumerate}
\item[(1)] Let $w\in W$ and 
let $w=s_1\cdots s_l$ be a representation.
If $\ell(w)<l$, then $w=s_1\cdots \hat{s_i} \cdots \hat{s_j} \cdots s_l$ 
for some $1\le i<j\le l$. 
\item[(2)] For each $w \in W$ and $s \in S$, 
$\ell(ws)$ equals either $\ell(w)+1$ or $\ell(w)-1$, and
$\ell(sw)$ also equals either $\ell(w)+1$ or $\ell(w)-1$.
\item[(3)] For each $w\in W$, 
$S(w)$ is a spherical subset of $S$, i.e., 
$W_{S(w)}$ is finite.
\item[(4)] For each $w\in W$ and each spherical subset $T$ of $S$, 
there exists a unique element of longest length in $W_Tw$. 
Here $v\in W_Tw$ is the element of longest length in $W_Tw$ 
if and only if $\ell(tv)<\ell(v)$ for any $t\in T$.
Then $\ell(v)=\ell(vw^{-1})+\ell(w)$.
\end{enumerate}
\end{Lemma}

\begin{Lemma}\label{lem3}
Let $(W,S)$ be a Coxeter system, 
let $w\in W$, let $s\in S$ such that $\ell(ws)=\ell(w)+1$, 
and let $T$ be a spherical subset of $S$.
By Lemma~\ref{lem1}~(4), 
there exist unique elements $x,x'\in W_T$ 
such that $xw$ and $x'ws$ are 
the elements of longest length in $W_Tw$ and $W_Tws$ respectively.
If $x=t_1\dots t_m$ is a reduced representation, 
then either $x'=x$ or 
$x'=t_1\cdots \hat{t_i} \cdots t_m$ for some $i\in\{1,\dots,m\}$.
Hence $\ell(x')\le\ell(x)$.
\end{Lemma}

\begin{proof}
By Lemma~\ref{lem1}~(2), 
either $\ell(xws)=\ell(xw)+1$ or $\ell(xws)=\ell(xw)-1$.

We first suppose that $\ell(xws)=\ell(xw)+1$.
Since $xw$ is the element of longest length in $W_Tw$, 
$\ell(txw)<\ell(xw)$ for any $t\in T$ by Lemma~\ref{lem1}~(4).
Then for each $t\in T$,
$$ \ell(txws)\le \ell(txw)+1<\ell(xw)+1=\ell(xws).$$
Hence $\ell(txws)<\ell(xws)$ for any $t\in T$.
Thus $xws$ is the element of longest length in $W_Tws$, 
i.e., $x'=x$.

Next we suppose that $\ell(xws)=\ell(xw)-1$.
Since $\ell(ws)=\ell(w)+1$ and $\ell(xw)=\ell(x)+\ell(w)$, 
$$ xws=(t_1\cdots t_m)ws=(t_1\cdots \hat{t_i} \cdots t_m)w $$
for some $i\in\{1,\dots,m\}$ by Lemma~\ref{lem1}~(1).
Now $xw$ is the element of longest length in $W_Tw$.
Here 
$$ W_Tw=W_T(t_1\cdots \hat{t_i} \cdots t_m)w=W_Txws=W_Tws.$$
Hence $xw$ is the element of longest length in $W_Tws$.
Since $xw=(t_1\cdots \hat{t_i} \cdots t_m)ws$, 
we obtain $x'=t_1\cdots \hat{t_i} \cdots t_m$.
\end{proof}

\begin{Definition}
Let $(W,S)$ be a Coxeter system.
For each $w \in W$, 
we define $S(w)= \{s \in S \,|\, \ell(ws) < \ell(w)\}$.
For a subset $T \subset S$, 
we also define $W^T= \{w \in W \,|\, S(w)=T \}$. 
\end{Definition}

The following lemma was proved in \cite{Ho2}.

\begin{Lemma}[{\cite[Lemma~2.5]{Ho2}}]\label{lem4}
Let $(W,S)$ be a Coxeter system, $w\in W$ and $s_0\in S$. 
Suppose that $m(s_0,t)\ge 3$ for each $t\in S(w)$ and 
that $m(s_0,t_0)=\infty$ for some $t_0\in S(w)$.
Then $ws_0\in W^{\{s_0\}}$.
\end{Lemma}

We can obtain the following lemma by 
the same argument as the proof of \cite[Lemma~4.2]{Ho}.

\begin{Lemma}\label{lem5}
Let $(W,S)$ be a Coxeter system 
and let $\alpha\in\partial\Sigma(W,S)$.
Then there exists a sequense $\{s_i\}\subset S$
such that $s_1\cdots s_i$ is reduced and
$d(s_1\cdots s_i,\Image\xi_{\alpha})\le N$ for each $i\in\N$, 
where 
$N$ is the diameter of $K(W,S)$ in $\Sigma(W,S)$ and 
$\xi_{\alpha}$ is the geodesic ray in $\Sigma(W,S)$ such that 
$\xi_{\alpha}(0)=1$ and $\xi_{\alpha}(\infty)=\alpha$.
\end{Lemma}

\section{Proof of the main theorem}

Using the lemmas above, we prove Theorem~\ref{Thm}.

\begin{proof}[Proof of Theorem~\ref{Thm}]
Let $(W,S)$ be a Coxeter system.
Suppose that there exist a maximal spherical subset $T$ of $S$ 
and an element $s_0\in S$ such that $m(s_0,t)\ge 3$ for each $t\in T$
and $m(s_0,t_0)=\infty$ for some $t_0\in T$.
Let $\alpha\in\partial\Sigma(W,S)$.
By Lemma~\ref{lem5}, 
there exists a sequense $\{s_i\}\subset S$
such that $s_1\cdots s_i$ is reduced and
$d(s_1\cdots s_i,\Image\xi_{\alpha})\le N$ for each $i\in\N$, 
where 
$N$ is the diameter of $K(W,S)$ in $\Sigma(W,S)$ and 
$\xi_{\alpha}$ is the geodesic ray such that 
$\xi_{\alpha}(0)=1$ and $\xi_{\alpha}(\infty)=\alpha$.
Let $w_i=s_1\cdots s_i$.
For each $i$, 
there exists a unique element $x_i\in W_T$ 
such that $x_iw_i$ is the element of longest length in $W_Tw_i$
by Lemma~\ref{lem1}~(4).
Now $w_{i+1}=w_is_{i+1}$ and $\ell(w_{i+1})=\ell(w_i)+1$.
By Lemma~\ref{lem3}, 
$\ell(x_{i+1})\le\ell(x_i)$ for every $i$.
Hence there exists a number $n$ 
such that $\ell(x_i)=\ell(x_{i+1})$ for each $i\ge n$.
Then $x_i=x_{i+1}$ for every $i\ge n$ by Lemma~\ref{lem3}.
Let $x=x_n$.
Then
$xw_i$ is the element of longest length in $W_Tw_i$ for each $i\ge n$.
Since $\ell(txw_i)<\ell(xw_i)$ for any $t\in T$, 
$T\subset S((xw_i)^{-1})$.
Here $S((xw_i)^{-1})$ is a spherical subset of $S$ by Lemma~\ref{lem1}~(3) 
and $T$ is a maximal spherical subset of $S$.
Hence $S((xw_i)^{-1})=T$ for each $i\ge n$.
By Lemma~\ref{lem4}, 
$(s_0xw_i)^{-1}\in W^{\{s_0\}}$ and 
$(t_0s_0xw_i)^{-1}\in W^{\{t_0\}}$ for each $i\ge n$.
Hence $(W^{\{t_0\}})^{-1}$ contains 
the sequence $\{t_0s_0xw_i\}_{i\ge n}$ 
which converges to $t_0s_0x\alpha$.
By the proof of \cite[Theorem~4.1]{Ho2}, 
$Wt_0s_0x\alpha$ is dense in $\partial\Sigma(W,S)$.
Here $Wt_0s_0x\alpha=W\alpha$.
Hence $W\alpha$ is a dense subset of $\partial\Sigma(W,S)$.
Thus every orbit $W\alpha$ is dense in $\partial\Sigma(W,S)$, 
that is, $\partial\Sigma(W,S)$ is minimal.
\end{proof}

%

%
\end{document}